\documentclass [12pt] {report}
\usepackage{amssymb}
\newcommand {\D}[2] {\displaystyle\frac{\partial{#1}}{\partial{#2}}}

\newcommand {\ga} {\gamma}
\newcommand {\la} {\lambda}

\newcommand {\de} {\delta}
\newcommand {\prtl} {\partial}
\newcommand {\fr} {\displaystyle\frac}

\newcommand {\be} {\begin{equation}}
\newcommand {\ee} {\end{equation}}
\newcommand {\ba} {\begin{array}}
\newcommand {\ea} {\end{array}}
\newcommand {\bp} {\begin{picture}}
\newcommand {\ep} {\end{picture}}
\newcommand {\bc} {\begin{center}}
\newcommand {\ec} {\end{center}}
\newcommand {\bt} {\begin{tabular}}
\newcommand {\et} {\end{tabular}}
\newcommand {\lf} {\left}
\newcommand {\rg} {\right}

\newcommand {\bls} {\baselineskip}

\newcommand {\cF} {{\cal F}}

\newcommand {\cR} {{\cal R}}

\newcommand {\ses} {\medskip}

\newcommand {\e} {\mathop{\rm e}\nolimits}

\newcommand {\nin} {\noindent}

\newcommand {\De} {\Delta}

\usepackage{epsfig, latexsym,graphicx}
 \usepackage{amsmath}

\def\2#1#2#3{{#1}_{#2}\hspace{0pt}^{#3}}
\def\3#1#2#3#4{{#1}_{#2}\hspace{0pt}^{#3}\hspace{0pt}_{#4}}
\newcounter{sctn}
\def\sec#1.#2\par{\setcounter{sctn}{#1}\setcounter{equation}{0}
                  \noindent{\bf\boldmath#1.#2}\bigskip\par}

\begin {document}

\begin {titlepage}

 \ses\ses\ses\ses

\vspace{0.1in}

\begin{center}

{\large \bf FINSLEROID--FINSLER  SPACE AND SPRAY

\ses\ses

 COEFFICIENTS}

\end{center}

\vspace{0.3in}

\begin{center}

\vspace{.15in} {\large G.S. ASANOV\\} \vspace{.25in}
{\it Division of Theoretical Physics, Moscow State University\\
119992 Moscow, Russia\\
{\rm (}e-mail: asanov@newmail.ru{\rm )}} \vspace{.05in}

\end{center}

\begin{abstract}

 In the previous work,
 the notion of the  Finsleroid--Finsler space  have been formulated and
the necessary and sufficient conditions for the
 space to be of the Landsberg type have been found.
In the present paper, starting with  particular spray coefficients, we demonstrate
how the Landsberg
condition can  explicitly appear in case of the Finsleroid--type metric function.
Calculations are supplementing by a convenient special Maple--program.
The general form of the associated geodesic spray coefficients is presented for such metric function
under the condition
of constancy of the Finsleroid charge.

\ses
\nin {\it  Key words}: Finsler geometry, metric spaces, spray.

\end{abstract}

\end{titlepage}

\vskip 1cm


\vskip 1cm

\setcounter{sctn}{1}
\setcounter{equation}{0}

 \bc
{\large 1. Introduction }
 \ec

\bigskip

Continuing the previous work [1-3] dealt with the Finsleroid--Finsler spaces,
we below clarify how the spray notion may entail the Landsberg--type Finsler space.

A {\it spray} ${\bf G}$ on an $N$--dimensional smooth
manifold $M$ is a smooth vector field on the slit tangent bundle $TM\backslash{0}$
 expressed in terms of a standard local coordinate system $(x^i,y^i)$ in $TM$ according to the representation
\be
{\bf G}=y^i\D{}{x^i}-G^i(y)\D{}{y^i}.
\ee
Spray spaces are generalized vector spaces which deep meaning is underlined by the property
that a spray ${\bf G}$ on $M$ determines a collection of geodesics in $M$, according to the differential
equation
\be
\fr{d^2c^i}{dt^2}+G^i\Bigl(\fr{dc}{dt}\Bigr)=0
\ee
for curves $c~:(a,b)\to M$ parametrized by $t$.
The theory of sprays bears close relation on the path spaces.
An interesting $S$-curvature can be associated with the spray concept.
The {\it Finsler geodesic spray}
 is the notion which is the adaptation of the general spray notion to the structure of Finsler spaces,
by using the Finslerian Christoffel symbols $\ga^k{}_{ij}$
and prescribing the equality
\be
G^k_{\{\text{Finsler}\}}=\ga^k{}_{ij}y^iy^j.
\ee
On the basis of  these coefficients the Finsler connection and curvature can consistently
be constructed by following known
methods (see [4-6]).

Suppose we are given on $M$ a Riemannian metric ${\cal S}=S(x,y)$ and a 1-form  $b=b(x,y)$ of the unit Riemannian length.
With respect to   local coordinates $x^i$ in the Riemannian space  $\cR_N=(M,{\cal S})$
we have the local representations $ b=b_i(x)y^i,$
 $S= \sqrt{a_{ij}(x)y^iy^j},$
 and
 \be
||b||_x~:=\sqrt{a^{ij}(x)b_ib_j}=1
\ee
with the tensor $a^{ij}$ reciprocal to the input
$a_{ij}$.
We shall construct from the covariant vector $b_i$ the contravariant vector $b^i$
 according to the Riemannian rule $b^i=a^{ij}b_j$.

\ses


We
also  introduce the tensor
\be
r_{ij}(x)~:=a_{ij}(x)-b_i(x)b_j(x),
\ee
obtaining the decomposition
\be
S^2=b^2+q^2
\ee
in terms of the scalar
\be
q~:=\sqrt{r_{ij}(x)y^iy^j}.
\ee
In many cases it is convenient to use the variables
\be
u_i~:=a_{ij}y^j,
\ee
\ses
\be
v^i~:=y^i-bb^i, \qquad v_m~:=u_m-bb_m=r_{mn}y^n\equiv r_{mn}v^n\equiv a_{mn}v^n,
\ee
 Notice that
\be
r^i{}_n~:=a^{im}r_{mn}=\de^i{}_n-b^ib_n=\D{v^i}{y^n}
\ee
(the $\de^i{}_n$ stands for the Kronecker symbol),
\be
v_ib^i=v^ib_i=0, \qquad r_{ij}b^j=r^i{}_jb^j=b_ir^i{}_j=0,
\ee
\ses
\be
u_iv^i=v_iy^i=q^2,
\ee
and
\be
\D b{y^i}=b_i, \qquad \D q{y^i}=\fr{v_i}q.
\ee
We  comply with the notation adopted in [2].

Under these conditions,
it seems attractive to take three scalars $c_1(x), c_2(x), c_3(x)$
and propose to consider on  the Riemannian space
 $\cR_N=(M,{\cal S})$  the spray given by the coefficients
\be
G^i=
c_1(x)\fr1q y^jy^h(\nabla_jb_h)v^i
+c_2(x)y^hb^j(\nabla_jb_h)v^i
+c_3(x)qf^i+a^i{}_{km}y^ky^m,
\ee
which are  such that the difference
$
G^i-a^i{}_{km}y^ky^m
$
involves all the  crucial terms linear in
 $ \nabla_jb_h.
 $
Here,  the nabla  means the covariant derivative in terms of the Riemannian space
 $\cR_N=(M,{\cal S})$;
 $ a^i{}_{nm} $ stands for
the  Riemannian  Christoffel symbols constructed from the tensor $a_{ij}(x)$;
the  notation
\ses
\be
f^i=f^i{}_ny^n,\qquad
f^i{}_n=a^{ik}f_{kn}, \qquad
f_{mn}=
\nabla_mb_n-\nabla_nb_m
\equiv \D{ b_n}{x^m}-\D {b_m}{x^n}
\ee
is used.


In Section 2 we consider  particular $G^i$ which reveal
the astonishing property of  the nullification (2.11)--(2.12) for contractions.
We call them  the  Landsberg--type spray coefficients, because, under the structural condition
(2.13),  Finsler metric functions inducing such
coefficients must produce the Landsberg--type spaces.

In Section 3 we demonstrate how the use of the generating Finsleroid--Finsler metric functions
 induces explicitly such coefficients  $G^i$.


In Section 4
the general representation of  the geodesic spray coefficients in case of the Finsleroid--Finsler space
with $g=const$ is given.
The representation  is obviously the kernel from which all the significant spray implications in such spaces are to be
grown up.

In Section 5 we present the Maple--program which verifies the   Landsberg--type spray coefficients.

The paper ends by Appendix A  in which the basic formulas and definitions of
the Finsleroid--Finsler space are summarized up.

Key  Propositions 1,2,3 are motivated and proven.


\setcounter{sctn}{2}
\setcounter{equation}{0}

\bc
{\large  2.  Landsberg--type spray coefficients }
\ec

\ses\ses

Under the condition
\be
\nabla_jb_i
=k(a_{ij}-b_ib_j)
\ee
with $ k=k(x),$
the coefficients (1.14)
reduce to
\be
G^i=cq(y^i-bb^i)+a^i{}_{km}y^ky^m\equiv gqkv^i+a^i{}_{km}y^ky^m,
\ee
where $c=c(x)$ is the  scalar that is obtained by
\be
c=c_1k.
\ee
From (2.2) we can readily calculate the entailed coefficients
\be
 G^i{}_k:~=\D{ G^i}{y^k}, \qquad  G^i{}_{km}:~=\D{ G^i{}_k}{y^m},
\qquad  G^i{}_{kmn}:~=\D{ G^i{}_{km}}{y^n}
\ee
by applying the rules (1.13), obtaining the representations
\be
 G^i{}_k=
\fr cq
\biggl[
(y^i-bb^i)(u_k-bb_k)+q^2(\de_k{}^i-b_kb^i)
\biggr]+2a^i{}_{km}y^m,
\ee
\ses\ses
$$
 G^i{}_{km}=\fr {c}q\biggl[
(a_{km}-b_kb_m)(y^i-bb^i)
-
\fr 1{q^2}(y^i-bb^i)(u_k-bb_k)
(u_m-bb_m)
$$
\ses
\be
+
(u_m-bb_m)(\de_k{}^i-b_kb^i)
+(u_k-bb_k)(\de_m{}^i-b_mb^i)
\biggr]
+2a^i{}_{km},
\ee
and
$$
G^i{}_{kmn}=
-\fr {c}{q^3}(u_n-bb_n)\Biggl[
(a_{km}-b_kb_m)v^i
-
\fr 1{q^2}v^i(u_k-bb_k)
(u_m-bb_m)
$$
\ses
$$
+
(u_m-bb_m)(\de_k{}^i-b_kb^i)
+(u_k-bb_k)(\de_m{}^i-b_mb^i)
\Biggr]
$$
\ses
$$
+\fr {c}q\Biggl[
(a_{km}-b_kb_m)(\de_n{}^i-b_nb^i)
+
\fr 2{q^4}(u_n-bb_n)v^i(u_k-bb_k)
(u_m-bb_m)
$$
\ses
$$
-
\fr 1{q^2}\Bigl((\de_n{}^i-b_nb^i)(u_k-bb_k)(u_m-bb_m)+v^i(a_{kn}-b_kb_n)(u_m-bb_m)+
v^i(u_k-bb_k)(a_{mn}-b_mb_n)
\Bigr)
$$
\ses
$$
+
(a_{mn}-b_mb_n)(\de_k{}^i-b_kb^i)
+(a_{kn}-b_kb_n)(\de_m{}^i-b_mb^i)
\Biggr].
$$


The explicitly symmetric form of the latter coefficients reads
$$
G^i{}_{kmn}=
\fr {3c}{q^5}
v^iv_kv_mv_n
-\fr {c}{q^3}\Biggl[
(\de_k{}^i-b_kb^i)v_mv_n
+
(\de_m{}^i-b_mb^i)v_kv_n
+
(\de_n{}^i-b_nb^i)v_kv_m
$$
\ses
$$
+v^i\Bigl(
(a_{km}-b_kb_m)v_n+(a_{kn}-b_kb_n)v_m+(a_{mn}-b_mb_n)v_k
\Bigr)
\Biggl]
$$
\ses
\be+\fr {c}{q}
\Biggl[
(\de_k{}^i-b_kb^i)(a_{mn}-b_mb_n)
+
(\de_m{}^i-b_mb^i)(a_{kn}-b_kb_n)
+(\de_n{}^i-b_nb^i)(a_{km}-b_kb_m)
\Bigr)
\Biggl],
\ee
or
\ses\\
\be
G^i{}_{kmn}=
\fr{c}{q}(\eta_k{}^i\eta_{mn}+\eta_m{}^i\eta_{kn}+\eta_n{}^i\eta_{km}),
\ee
where the $\eta$--tensors are given by
\be
\eta^i{}_j~:=r^i{}_j-\fr1{q^2}v^iv_j, \qquad \eta_{ij}~:=r_{ij}-\fr1{q^2}v_iv_j\equiv a_{in}\eta^n{}_j
\ee
(the formulas (1.7)--(1.10) have been used).
Because of the nullifications
\be
b_i\eta^i{}_j=u_i\eta^i{}_j=0
\ee
(see the formulas (1.11)--(1.12))
the obtained coefficients (2.8) fulfill the identities
\be
b_i G^i{}_{kmn}=0
\ee
and
\be
 u_iG^i{}_{kmn}=0.
\ee


Suppose a  Finslerian metric function $F(x,y)$ be obtainable  from a function $\breve F(S,b)$ of the 1-form $b$ and
the Riemannian metric function $S$, such that
\be
F(x,y)=\breve F(S(x,y), b(x,y)).
\ee
Then it is obvious that  covariant vectors $\{y_i\}$ produced by the function $F$ according to the conventional
Finsler rule $y_i=\fr12\partial F^2(x,y)/\partial y^i$ are linear combinations of $b_i$ and $u_i$,  that is,
the equality
\be
y_i=p_1b_i+p_2u_i
\ee
holds with two scalars
$p_1,p_2$.
Noting the vanishings  (2.11) and (2.12),
 we are justified to claim the following.

\ses\ses

\nin
{\bf Proposition 1.} {\it Suppose a Finsler metric function $F(x,y)$ entail the spray coefficients of the
form \rm(2.2). \it If also the function $F$  is of the structure \rm(2.13),
\it
then the function $F$ produces the identity
\be
y_iG^i{}_{kmn}=0
\ee
and, hence, a Landsberg--case Finsler space.
}

\ses\ses

Because of this observation, we introduce the following.

\ses

{\bf Definition.} The coefficients $G^i$ given by  the representation (2.2)
are called the {\it  Landsberg--type spray coefficients}.


{\it In the two--dimensional case,}
\be
N=2,
\ee
such a unit 1-form  $e=e_j(x)y^j$ exists
that
\be
a_{ij}=e_ie_j+b_ib_j,
\ee
whence
the definitions (1.5), (1.7), and (1.9)  reduce to
\be
r_{ij}=e_ie_j,
\ee
\be
 q=|e|,
 \ee
 \be
 v^i=ee^i,
\ee
the $\eta$--tensors (2.9) vanish
\be
\eta_k{}^i=\eta_{kn}=0,
\ee
and
the implication
\be
(N=2)\quad \to \quad G^i{}_{kmn}=0
\ee
is applicable to (2.8) {\it  independently of the value of the scalar $c$.}
Therefore, in the dimension $N=2$ the coefficients (2.2), (2.5), and (2.6) reduce to
\be
G^i=ce|e|e^i+a^i{}_{km}y^ky^m,
\ee
\ses
\be
 G^i{}_k=
2c|e|e^ie_k
+2a^i{}_{km}y^m,
\ee
and
\be
 G^i{}_{km}=\fr{2ce}{|e|}e^ie_ke_m+2a^i{}_{km}.
\ee
 The latter coefficients are {\it independent} of vectors $y$,
 thereby corresponding to the Berwald case.

The formula (2.2) obtained for the spray coefficients,
 as well as the very condition (2.1),
 is applicable {\it in any} dimension
$N\ge2$. The right--hand side of the formula involves $q$ which is a square root of a quadratic form of rank
$N-1$, so that
in the dimensions $N\ge 3$ the coefficients can not be quadratic in vectors $y$ (unless the Riemannian case occurs),
--- this note
may be regarded as  the reason proper why the Landsberg case treated does not degenerate to the Berwald case
at $N\ge3$.
In the dimension $N=2$,
however,
the quadratic form mentioned is a square of the 1-form $e$ introduced above,
hence the square root is extracted up (see (2.19)), leaving us
with the expression (2.23) quadratic in vectors $y$, that is  with  the Berwald case.
In the Finsleroid--Finsler space of the dimension $N=2$
the associated main scalar $I$ proves to be $ I= I(x)=|g(x)|$ (cf. p. 26 in [2]).




\setcounter{sctn}{3}
\setcounter{equation}{0}

\bc
{\large 3. Use of  generating  metric functions in the Finsleroid--Finsler case}
\ec

\ses\ses

Let us  inquire into whether the Finsleroid--Finsler metric function $K$
with $g=const$ may fulfill the conditions which underlined Proposition 1.

Accordingly, we use a constant $g$ ranging over
$
-2<g<2,
$
 together with the notation
\be
h=\sqrt{1-\fr14g^2}, \qquad
G=g/h.
\ee
The respective Finsleroid--Finsler metric function $K$ does belong to the class (2.13).

Put
$
w=q/b
$
whenever $b\ne0$
and  rewrite the function $K$ in the form
\be
K=bV(w),
\ee
where the {\it generating  metric function} $V(w)$ is smooth of  the class $C^{\infty}$ on all the region
\be
 I_w= (-\infty,\infty).
 \ee
In terms of the quadratic form
\be
Q(w)=1+gw+w^2
\ee
we have
\be
V(w)=\sqrt{Q(w)}\,\e^{\frac12G\Phi(w)},
\ee
where
\be
\Phi(w)=
\fr{\pi}2+\arctan \fr G2-\arctan\Bigl(\fr{w+\fr g2}{h}\Bigr),
\qquad  {\rm if}
 \quad b\ge 0,
\ee
and
\be
\Phi(w)=
-\fr{\pi}2+\arctan \fr G2-\arctan\Bigl(\fr{w+\fr g2}{h}\Bigr),
\qquad  {\rm if}
\quad
0\ge b,
\ee
We obtain
\be
V'=wV/Q,
\qquad
V''=V/Q^2,
\ee
\ses
\be
(V^2/Q)'=-gV^2/Q^2,  \qquad (V^2/Q^2)'=-2(g+w)V^2/Q^3,
\qquad \Phi'=-h/Q
\ee
\bigskip\\
and also
\be
\fr12(V^2)'=wV^2/Q,
\qquad\quad
\fr12(V^2)''=(Q-gw)V^2/Q^2,
\qquad \fr14(V^2)'''=-gV^2/Q^3,
\ee
where the prime ($'$) denotes the differentiation with respect to~$w$.


If, alternatively,  we use the variable
$
s=b/S
$
and consider the function $K$ to read
\be
K=S\phi(s),
\ee
we obtain the {\it generating  metric function} $\phi(s)$ which is  smooth of the  class $C^{\infty}$ on the interval
\be
I_s= (-1,1),
\ee
with
\be
\phi(s)=
\sqrt{1+gs\sqrt{1-s^2}}\,\e^{\frac12G\Phi(s)},
\ee
where
\be
\Phi(s)=
\fr{\pi}2+\arctan \fr G2-\arctan\Bigl(\fr{\sqrt{1-s^2}+\fr g2s}{hs}\Bigr),
\qquad  {\rm if}
 \quad 1>s\ge 0,
\ee
and
\be
\Phi(s)=
-\fr{\pi}2+\arctan \fr G2-\arctan\Bigl(\fr{\sqrt{1-s^2}+\fr g2s}{hs}\Bigr),
\qquad  {\rm if}
\quad
0\ge s>-1,
\ee
The limits at $s=0$ from the left and from the right are the same value
$$ \phi(0)=\e^{\frac12G\arctan (G/2)}. $$
Evaluating derivatives yields merely
\ses\\
\be
\phi'=\fr{g\sqrt{1-s^2}\,\e^{\frac12G\Phi(s)}}
{\sqrt{1+gs\sqrt{1-s^2}}}, \qquad
\phi''=-\fr{gs\,\e^{\frac12G\Phi(s)}}{\sqrt{1-s^2}
\Bigl(\sqrt{1+gs\sqrt{1-s^2}}\Bigr)^3},
\ee
\ses
\be
\phi(\phi-s\phi')=\e^{G\Phi(s)},
 \qquad \phi-s\phi'+(1-s^2)\phi''=\fr{\e^{\frac12G\Phi(s)}}
{\Bigl(\sqrt{1+gs\sqrt{1-s^2}}\Bigr)^3},
\ee
and
\ses\\
\be
\fr{\phi\phi'-s(\phi\phi''+\phi'\phi')}
{\phi\Bigl((\phi-s\phi')+(1-s^2)\phi''\Bigr)}
=\fr g{\sqrt{1-s^2}},
\qquad
\fr{\phi''}{(\phi-s\phi')+(1-s^2)\phi''}
=-\fr {gs}{\sqrt{1-s^2}},
\ee
\ses
\ses
\be
\fr{(\phi-s\phi')^2} {\phi[\phi-s\phi'+(1-s^2)\phi'']}=1,
\ee
\ses\\
where the prime ($'$) denotes the differentiation with respect to the variable $s$.
Singularities appear when $|s|\to 1$.


\ses

Using each of the two generating metric functions, $V(w)$ or $\phi(s)$, it is easy to observe that
with the symmetry assumption
\be
\nabla_jb_i=\nabla_ib_j
\ee
the Finsler spray  coefficients
$
 G^i =
\ga^i{}_{nm}y^ny^m
$
(see the formula (1.3)) prove to be of the explicit form
\ses\\
\be
G^i=\fr{gy^my^n\nabla_nb_m}{\sqrt{S^2-b^2}}(y^i-b b^i)+a^i{}_{nm}y^ny^m.
\ee
If we plug here the condition (2.1),
we obtain the spray coefficients
\ses\\
\be
G^i=gk\sqrt{S^2-b^2}(y^i-b b^i)+a^i{}_{nm}y^ny^m
\ee
which  are   tantamount to the Landsberg--type
spray coefficients (2.2).

\ses\ses

This way  we have arrived at the following.

\ses\ses

{\bf Proposition 2}.
{\it Under the  assumption \rm(2.1), \it the Finsleroid--Finsler metric function $K$ with $g=const$
induces the Landsberg--type spray coefficients. The entailed coefficients
$
G^i{}_{kmn}
$
are of the simple form that is given by
}
(2.8),
 {\it with}
\be
c=kg.
\ee

\ses\ses

The Berwald case  corresponds to $k=0$ in dimensions $N\ge 3$, and holds uniquely in dimension $N=2$.

\ses\ses\ses\ses


\setcounter{sctn}{4}
\setcounter{equation}{0}

\bc
{\large 4. General form of  Finsleroid--Finsler geodesic spray coefficients}
\ec

\ses\ses

Straightforward calculations of the
 Finsleroid--Finsler Christoffel symbols
$\ga_{ikj}$ results in the following representation:
\be
2\ga_{ikj}=
\fr{\partial g_{kj}}{\partial g}\fr{\partial g}{\partial x^i}
+\fr{\partial g_{ik}}{\partial g}\fr{\partial g}{\partial x^j}
-\fr{\partial g_{ij}}{\partial g}\fr{\partial g}{\partial x^k}
$$
\ses
$$
+\fr{2gb^2}{Bq}(c_ig_{kj}+c_jg_{ik}-c_kg_{ij})
+\fr g{B}(q-\fr{b^2}q)(c_ia_{kj}+c_ja_{ik}-c_ka_{ij})\fr{K^2}B
$$
\ses\ses
$$
-\fr {g\fr{K^2}B}{B}
\Biggl[
c_i\biggl[(2gb+\fr{{S^2}}q)b_kb_j+
\fr1{q^2}\fr {S^2}q(b^2b_kb_j-bb_ku_j-bb_ju_k+u_ku_j)
\biggr]
$$
\ses\ses
$$
+\lf(q+\fr{b^2}q\rg)\Bigl(b_{k,i}(bb_j-u_j)+b_{j,i}(bb_k-u_k)\Bigr)
-gq^2(b_{k,i}b_j+b_{j,i}b_k)
\Biggr]
$$
\ses\ses
$$
-\fr {g\fr{K^2}B}{B}
\Biggl[
c_j\biggl[(2gb+\fr{{S^2}}q)b_kb_i+
\fr1{q^2}\fr {S^2}q(b^2b_kb_i-bb_ku_i-bb_iu_k+u_ku_i)
\biggr]
$$
\ses\ses
$$
+\lf(q+\fr{b^2}q\rg)\Bigl(b_{k,j}(bb_i-u_i)+b_{i,j}(bb_k-u_k)\Bigr)
-gq^2(b_{k,j}b_i+b_{i,j}b_k)
\Biggr]
$$
\ses\ses
$$
+\fr {g\fr{K^2}B}{B}
\Biggl[
c_k\biggl[(2gb+\fr{{S^2}}q)b_ib_j+
\fr1{q^2}\fr {S^2}q(b^2b_ib_j-bb_iu_j-bb_ju_i+u_iu_j)
\biggr]
$$
\ses\ses
$$
+\lf(q+\fr{b^2}q\rg)\Bigl(b_{i,k}(bb_j-u_j)+b_{j,k}(bb_i-u_i)\Bigr)
-gq^2(b_{i,k}b_j+b_{j,k}b_i)
\Biggr]
+\De.
\ee
Here, $S^2=b^2+q^2$, $u_i=a_{ij}y^j$, $b_{j,k}=\partial b_j/\partial x^k$,
 $c_i=y^kb_{k,i}$, and $\De$ symbolizes the summary of the terms which involve partial derivatives
of the input Riemannian metric tensor $a_{ij}$ with respect to the coordinate variables $x^k$.


Using this result, we obtain after due calculations the following.

\ses\ses

\nin
{\bf Proposition 3}.
{\it In the Finsleroid--Finsler  space under the only condition that
the  Finsleroid charge is a constant,  $ g=const$,
the induced spray coefficients $
 G^i =
\ga^i{}_{nm}y^ny^m
$
can  explicitly be written in the form \rm(1.14) \it with}
\be
 c_1=g, \qquad c_2=g^2,\qquad  c_3=-g,
\ee
so that
\be
G^i=
 g
\Bigl(\fr1q y^jy^h\nabla_jb_h+gy^hb^j\nabla_jb_h\Bigr)
v^i
-gqf^i+a^i{}_{km}y^ky^m.
\ee

\ses\ses

We have here $f^i=0$ if  the symmetry (3.20)
is assumed (see (1.15)).
If the condition (2.1) is plugged in (4.3), the spray coefficients (2.2) appear with $c$ given by (3.23).

It is remarkable to note that  the Finsleroid--Finsler metric function $K$ does not enter the right--hand side of (4.3).
The presence of the constant $g$ in the right--hand side of (4.3)  is the only trace of the function $K$ in the
spray coefficients $G^i$ obtained.


\setcounter{sctn}{5}
\setcounter{equation}{0}

\bc
{\large  5.  Maple--verification }
\ec

\ses\ses

Below we  check the vanishing $\dot{A}_{jkl}=0$  by the resource of the Maple10, using the formulas
\be
 G^i=cq(y^i-bb^i)+a^i_{km}y^ky^m,\qquad
\dot{A}_{jkl}=-\frac14y_i\frac{\partial^3 G^i}{\partial y^j\partial y^k\partial y^l}
\ee
with $ q=\sqrt{r_{ij}y^iy^j}$ and  $b=b_iy^i;$  $c$ is independent of $y$.
In the program, $ G^i$ will be denoted by gammas[i], $\dot{A}_{jkl}$ by dotA[j,k,l], and $b$ by bs.
\begin{verbatim}
> restart:
> N:=2:q:=sqrt(add(add(r[i,j]*y[i]*y[j],j=1..N),i=1..N)):
  bs:=add(b[i]*y[i],i=1..N):
  for i from 1 to N do
  c*q*(y[i]-bs*b[i])+add(add(a[i,j,k]*y[j]*y[k],j=1..N),k=1..N):
  gammas[i]:=eval(%,{seq(seq(r[i,j]=r[j,i],i=1..j),j=1..N)});
  end do:
\end{verbatim}

\nin
Apply $y_i=p_1b_i+p_2r_{ij}y^j$ (Eq. (2.14)) with arbitrary $p_1$ and $p_2$.

\begin{verbatim}
> for j from 1 to N do for k from 1 to N do for l from 1 to N do
  dotA[j,k,l]:=-1/4*factor(add( (add(r[i,a]*y[a]*p2,a=1..N)+b[i]*p1)
  *diff(diff(diff(gammas[i],y[l]),y[k]),y[j]),i=1..N));
  end do:end do:end do:
\end{verbatim}

\nin
Plug the symmetry $r_{ij}=r_{ji}$ and
simplify the arisen quantities  $\dot A_{jkl}$
by the use of the constrains
$ b_ib^i=1$ and $r_{ij}b^j=0.$

\begin{verbatim}
> for a1 from 1 to N do for a2 from 1 to N do for a3 from 1 to N do
  dotA[a1,a2,a3]:factor(eval(%,{seq(seq(r[i,j]=r[j,i],i=1..j),j=1..N)})):
  algsubs(add(b[i]^2,i=1..N)=1,%):factor(%);
  for j from 1 to N do
  algsubs(eval(add(r[i,j]*b[i],i=1..N)=0,
  {seq(seq(r[i,j]=r[j,i],i=1..j),j=1..N)}),%);end do:
  simplify(%);print(%);end do:end do:end do:
\end{verbatim}

\nin
The result  of the simplification is just the succession of zeros:
\[0\]\[0\]\[0\]\[0\]\[0\]\[0\]\[0\]\[0\]
This result supports Proposition 1 of Section 2.


The calculation times at the dimensions $N=2$ and $N=3$ are short.



\setcounter{equation}{0}

\bc
{\large  Appendix A:  Involved Finsleroid--Finsler representations}
\ec

We introduce on the manifold $M$
a scalar  $g=g(x)$
 subject to ranging
\be
-2<g(x)<2,
\ee
and apply  the convenient notation
\be
h=\sqrt{1-\fr14g^2}, \qquad
G=g/h.
\ee
\ses

The {\it  characteristic
quadratic form}
\be
B(x,y) :=b^2+gqb+q^2
\equiv\fr12\Bigl[(b+g_+q)^2+(b+g_-q)^2\Bigr]>0
\ee
where $ g_+=\fr12g+h$ and $ g_-=\fr12g-h$,
is of the negative discriminant
\be
D_{\{B\}}=-4h^2<0
\ee
and, therefore, is positively definite.

\ses\ses

 {\large Definition}. The scalar function $K(x,y)$ given by the formulas
\be
K(x,y)=
\sqrt{B(x,y)}\,J(x,y)
\ee
and
\be
J(x,y)=\e^{\frac12G\Phi(x,y)},
\ee
where
\be
\Phi(x,y)=
\fr{\pi}2+\arctan \fr G2-\arctan\Bigl(\fr{L(x,y)}{hb}\Bigr),
\qquad  {\rm if}  \quad b\ge 0,
\ee
and
\be
 \Phi(x,y)= -\fr{\pi}2+\arctan \fr
G2-\arctan\Bigl(\fr{L(x,y)}{hb}\Bigr), \qquad  {\rm if}  \quad
b\le 0,
\ee
 with
 \be
 L(x,y) =q+\fr g2b,
\ee
\ses\\
is called
the {\it  Finsleroid--Finsler  metric function}.

\ses\ses

The positive (not absolute) homogeneity holds fine: $K(x,\la y)=\la K(x,y)$ for all $\la >0$.

In the limit $g\to 0$,
the definition  degenerates to the
 input Riemannian metric function:
\be
K|_{_{g=0}}=S.
\ee


\ses\ses

 {\large  Definition}.  The arisen  space
\be
\cF\cF^{PD}_g :=\{\cR_{N};\,b(x,y);\,g(x);\,K(x,y)\}
\ee
is called the
 {\it Finsleroid--Finsler space}.

\ses\ses

 {\large  Definition}. The space $\cR_N$ entering the above definition is called the {\it associated Riemannian space}.

\ses\ses

{\large Definition}.\, Within  each tangent space $T_xM$, the Finsleroid--metric function $K(x,y)$  produces the {\it Finsleroid}
 \be
 \cF^{PD}_{g\,\{x\}}:=\{y\in  \cF^{PD}_{g\,\{x\}}: y\in T_xM , K(x,y)\le 1\}.
  \ee

We calculate from the function $K$ the
 covariant tangent vector $\hat y=\{y_i\}$ and
the  Finslerian metric tensor $\{g_{ij}\}$, by making  use of the  conventional  Finslerian  rules
\be
y_i :=\fr12\D{K^2}{y^i}, \qquad
g_{ij} :
=
\fr12\,
\fr{\prtl^2K^2}{\prtl y^i\prtl y^j}
=\fr{\prtl y_i}{\prtl y^j},
\ee
obtaining
\be
y_i=(a_{ij}y^j+gqb_i)\fr{K^2}B
\ee
and
\be
g_{ij}=
\biggl[a_{ij}
+\fr g{B}\Bigl((gq^2-\fr{bS^2}q)b_ib_j-\fr bqu_iu_j+
\fr{ S^2}q(b_iu_j+b_ju_i)\Bigr)\biggr]\fr{K^2}B,
\ee
where  the notation (1.8) has been used.
The reciprocal components $(g^{ij})=(g_{ij})^{-1}$ read
\be
g^{ij}=
\biggl[a^{ij}+\fr gq(bb^ib^j-b^iy^j-b^jy^i)+\fr g{Bq}(b+gq)y^iy^j
\biggr]\fr B{K^2},
\ee


In terms of the variables (1.9) we obtain the representations
\be
y_i=\Bigl(v_i+(b+gq)b_i\Bigr)\fr{K^2}B,
\ee
\ses
\be
g_{ij}=
\biggl[a_{ij}
+\fr g{B}\Bigl (q(b+gq)b_ib_j+q(b_iv_j+b_jv_i)-b\fr{v_iv_j}q\Bigr)\biggr]\fr{K^2}B,
\ee
and
\be
g^{ij}=
\biggl[a^{ij}+\fr gB\Bigl(-bqb^ib^j-q(b^iv^j+b^jv^i)+(b+gq)\fr{v^iv^j}q\Bigr)
\biggr]\fr B{K^2}
\ee
which are alternative to (A.14)--(A.16).

The determinant of the metric tensor  is the smooth and positive function as follows:
\be
\det(g_{ij})=\bigl(\fr{K^2}B\bigr)^N\det(a_{ij})>0.
\ee
For the component of the contracted Cartan tensor
we find
\be
A_i=\fr {NK}2g\fr1{q}(b_i-\fr b{K^2}y_i)=\fr {NK}2g\fr1{qB}(q^2b_i- bv_i).
\ee

Since
\be
\fr{v^iv^j}q\to 0 \quad {\text{when}}\quad v^i\to 0
\ee
(notice the definition (1.7) of $q$) the components  $y_i$, $g_{ij}$ and $g^{ij}$, as given by (A.17),
(A.18), and (A.19),
are smooth on all the slit tangent bundle
$TM\backslash 0$.
However, the components (A.21) are singular at $v^i=0$.
 Therefore, on $TM\backslash 0$ the Finsleroid--Finsler space is smooth of the class $C^2$ and not of the class $C^3$,
 and at the same time the space
 is smooth of the class $C^{\infty}$
 on
$
TM\backslash\{0,-b,b\}.
$

We use  the
Riemannian Christoffel symbols
\be
a^k{}_{ij}~:=\fr12a^{kn}(\prtl_ja_{ni}+\prtl_ia_{nj}-\prtl_na_{ji})
\ee
($\prtl_j=\prtl/\prtl x^j$)
 given rise to by the associated Riemannian metric ${\cal S}$,
and also  the Finslerian Christoffel symbols
\be
\ga^k{}_{ij}~:=g^{kn}\ga_{inj}
\ee
with
\be
\ga_{inj}~:=\fr12(\prtl_jg_{ni}+\prtl_ig_{nj}-\prtl_ng_{ji}).
\ee


In terms of the  tensors
\be
{\cal H}_{mn}=\eta_{mn}\fr{K^2}B, \qquad{\cal H}_k{}^i=\eta_k{}^i
\ee
(cf. Eqs. (A.11) and (A.16)  in [2])
 the coefficients given by (2.8) take on the form
\be
G^i{}_{kmn}=
\fr{c}{q}({\cal H}_k{}^i{\cal H}_{mn}+{\cal H}_m{}^i{\cal H}_{kn}+{\cal H}_n{}^i{\cal H}_{km})\fr B{K^2},
\ee
and if we lower here the first index,  we obtain the {\it totally symmetric coefficients}
\be
G_{ikmn}=
\fr{c}{q}({\cal H}_{ik}{\cal H}_{mn}+{\cal H}_{im}{\cal H}_{kn}+{\cal H}_{in}{\cal H}_{km}).
\ee
The vanishings
\be
y^iG_{ikmn}=0, \qquad  b^iG_{ikmn}=0, \qquad  A^iG_{ikmn}=0
\ee
hold.


\ses
\ses

\end{document}